# ON THE VIBRATIONS OF LUMPED PARAMETER SYSTEMS GOVERNED BY DIFFERENTIAL-ALGEBRAIC EQUATIONS

S. DARBHA, K. B. NAKSHATRALA, AND K. R. RAJAGOPAL

ABSTRACT. In this paper, we consider the vibratory motions of lumped parameter systems wherein the components of the system cannot be described by constitutive expressions for the force in terms of appropriate kinematical quantities. Such physical systems reduce to a system of differential-algebraic equations, which invariably need to be solved numerically. To illustrate the issues with clarity, we consider a simple system in which the dashpot is assumed to contain a "Bingham" fluid for which one cannot describe the force in the dashpot as a function of the velocity. On the other hand, one can express the velocity as a function of the force.

## 1. INTRODUCTION

The traditional approach to obtaining the governing equations for the vibratory motions of a lumped parameter system of springs, masses and dashpots is to write down the balance of linear momentum for the system and to provide constitutive expressions for the forces in the springs and dashpots in terms of appropriate kinematical quantities. This leads to a differential equation for the motion which can be solved by a proper choice of the initial conditions for the system. However, in many instances, one is not in a position to provide the constitutive expressions for the forces in the components of the system in terms of the kinematical quantities, rather one finds that constitutive expressions can be provided for the appropriate kinematical quantity in terms of the forces. Two examples which immediately come to mind are a dashpot wherein the fluid behaves like a "Bingham fluid" or the dissipation due to Coulomb friction. In the case of a dashpot, wherein the fluid responds like a "Bingham fluid" (see Figure 1), we immediately observe that the force is not a function of the velocity, but the velocity is a function of the force. (If by a fluid, one means a body that cannot resist shear stress, then the notion of a "Bingham fluid" is untenable as there cannot be a yield stress below which the fluid resists shear stress.) It is also possible that the response is such that the relationship between the velocity and the force is not invertible (see Reference [13]). Similarly, it is possible that the energy storage mechanism is such that one cannot express the force in terms of the displacement but can only express the displacement in terms of the force. Such a situation arises when one comes across an energy storage mechanism described by the lumped parameter system depicted in Figure 2(a), which consists of a spring and an inextensible string in parallel. The relation between the force and displacement is that depicted in Figure 2(b). This response can be expressed as the displacement being a function of the force as given in Figure 2(c).

In general, the constituents of the lumped parameters may only be characterizable with implicit constitutive relations between the forces and kinematical quantities. Recently, Rajagopal [12] has discussed the general framework that arises when considering such lumped parameter systems. The problem reduces to the solution of a system of differential-algebraic equations. In this paper, we consider such a system of equations, for a rather simple lumped parameter system consisting of a mass, spring and dashpot, the fluid in the dashpot being characterized by the constitutive expressions for a "Bingham fluid", with a view towards explaining the interesting features of such a differential-algebraic system of equations. This study is just an initial foray into an area which holds much promise in view of the lumped parameter system being much more complicated. For instance, as we mentioned earlier, the components of the lumped parameter system may be such that one can only provide implicit constitutive relations to describe them, with neither the force nor the kinematical quantity being expressible in terms of the other. It is possible that a constituent of the lumped parameter might have both an energy storage and an energy dissipation mechanism built into them, that is the lumped object cannot be further split into an energy storing (spring) and an energy dissipating (dashpot) mechanism. That is, the object that is being "lumped" might be viscoelastic. An example of such a situation presents itself when the spring-dashpot system reflects the representation of a Maxwell fluid [10] or one that corresponds to more complicated rate type of fluids. The system could then be comprised of various such lumped objects and more than one kinematical quantity and force might be given in terms of implicit equations (see [12] for a more detailed discussion of such lumped parameter systems). Also complicating matters further could be dissipative systems that are Coulomb-like. Such systems, even within the framework of simpler components are far from being well understood (see [15, 1] for a discussion concerning the mathematical difficulties associated with the study of such systems).

In this short paper, we simplify the problem greatly in order to highlight the salient aspects of the problem under consideration.We study a simple mass, spring and dashpot system wherein the spring, to keep matters simple, is a linear spring and the dashpot is a "Bingham dashpot" whose response is given in Figure 4. We obtain the differential-algebraic system corresponding to the motion of such a system and solve it numerically for different initial conditions. We find interesting results that are in keeping with physical expectations. For instance, in the case of free vibrations, it is possible that the system could come to rest in a position wherein the spring is stretched, i.e., the equilibrium solution is not the solution corresponding to the unstretched spring. Thus, the equilibrium solution cannot be obtained by minimizing the "energy" associated with the system. This result is quite different from what one obtains in the case of a classical viscoelastic dashpot wherein the equilibrium state is the unstretched state of the spring. We also study the problem of the system being subject to different types of forced vibrations, the results being in keeping with physical expectations.

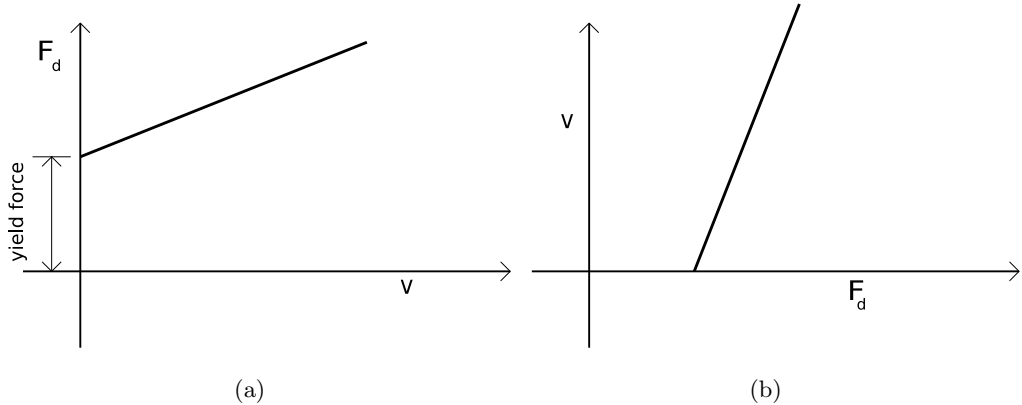

FIGURE 1. Bingham fluid model (a) force-velocity relationship (b) velocity as a *function* of the force

## 2. PROBLEM STATEMENT AND GOVERNING EQUATIONS

Consider a mass-spring-dashpot system as shown in Figure (3). We shall denote the deflection of the spring from its unstretched position by $x(t)$. We place an inertial frame of reference at the point occupied by the mass and hence, we may represent the displacement of the mass from this position by $x(t)$. Let $m$ denote the mass. Let $F_d(t)$ and $F_s(t)$, respectively, denote the forces applied by the dashpot and the spring on the mass, and let $F(t)$ denote the applied external force acting on the mass. An application of the balance of linear momentum yields the following equation:

$$\ddot{x}(t) = \frac{1}{m}\left(F(t) - F_d(t) - F_s(t)\right), \tag{1}$$

where a superposed dot denotes the time derivative. In the first-order form, the above equation takes the following form

$$\dot{v}(t) = \frac{1}{m}[F(t) - F_d(t) - F_s(t)], \tag{2a}$$

$$\dot{x}(t) = v(t), \tag{2b}$$

where $v(t)$ is the velocity of the mass.

One must now specify a constitutive relationship for the spring force $F_s(t)$ and dashpot force $F_d(t)$ in terms of the kinematical variables, usually $x(t)$ and $v(t)$, so that one can determine the

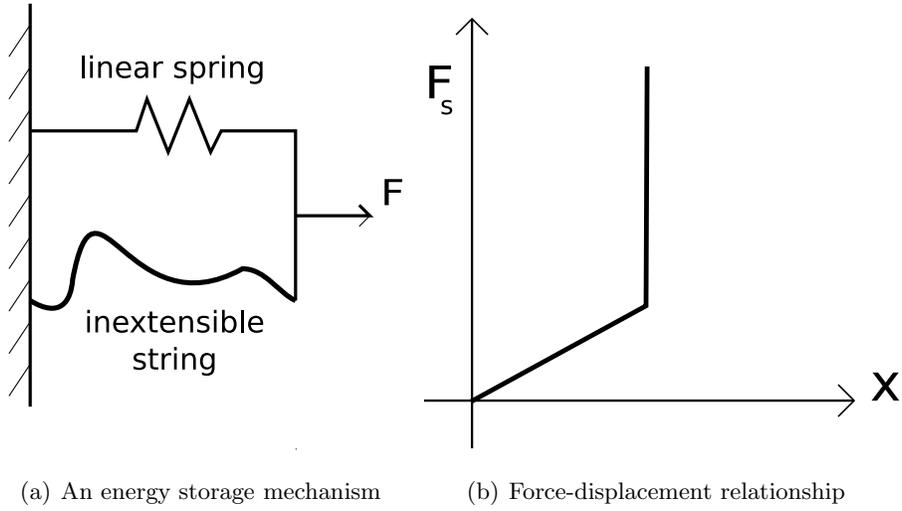

(a) An energy storage mechanism

(b) Force-displacement relationship

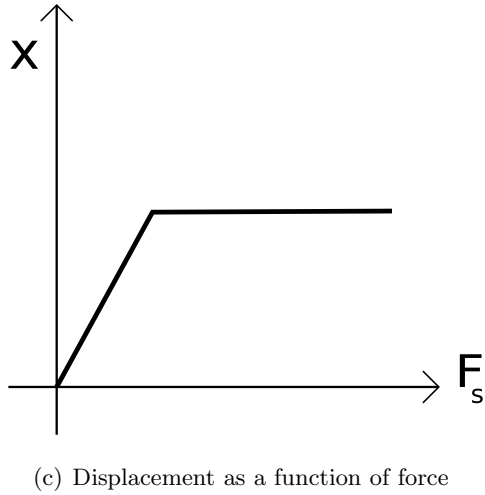

(c) Displacement as a function of force

FIGURE 2. A spring with an inextensible string in parallel

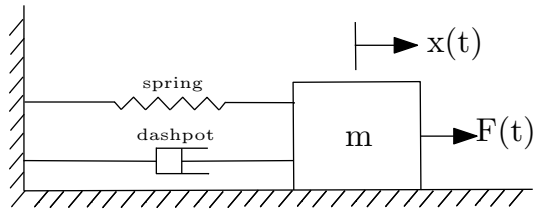

FIGURE 3. A mass-spring-dashpot system

response of the mass. Generally, these constitutive relationships usually take the form

(3a) $$0 = \alpha(x(t), F_s(t)),$$
(3b) $$0 = \beta(v(t), F_d(t)).$$

where $\alpha$ and $\beta$ are appropriate functions. However, one could have a much more complicated implicit constitutive relationship than (3a) and (3b). (For instance, the spring could be a viscoelastic solid though one invariably uses the term spring to denote a purely elastic body.) Then, the governing equations for the mass-spring-dashpot system may be written as

$$\dot{v}(t) = \frac{1}{m}[F(t) - F_s(t) - F_d(t)], \tag{4a}$$

$$\dot{x}(t) = v(t), \tag{4b}$$

$$0 = \alpha(x(t), F_s(t)), \tag{4c}$$

$$0 = \beta(v(t), F_d(t)). \tag{4d}$$

Now the variables are $x(t)$, $v(t)$, $F_d(t)$ and $F_s(t)$. In addition to the above equations, we need to specify (consistent) initial conditions. In the next subsection, we shall consider specific models for the spring and dashpot.

2.1. **Specific models for spring and dashpot.** In this paper we consider the case of spring, which is linear, and the dashpot is assumed to contain a classical Bingham fluid. The spring constant will be represented by $k$. The constitutive equation for the spring may be written as

$$x(t) = \frac{1}{k} F_s(t). \tag{5}$$

To provide the constitutive relationship for the dashpot containing a Bingham fluid, we require two *non-negative* constants: $\mu_s N_f$ and $\gamma$, where $\mu_s$ is coefficient of static friction, and $N_f$ is the normal force on the mass. The constitutive relationship (4d) takes the following explicit form:

$$v(t) = \begin{cases} 0 & |F_d(t)| \leq \mu_s N_f \\ \gamma\left(F_d(t) - \operatorname{sgn}[F_d(t)]\mu_s N_f\right) & |F_d(t)| > \mu_s N_f \end{cases} \tag{6}$$

where $\gamma > 0$ is the slope and $\mu_s N_f \geq 0$ is the threshold force (see Figure 4). For the case of $|F_d(t)| > \mu_s N_f$ it is easy to see that

$$\operatorname{sgn}[v(t)] = \operatorname{sgn}[F_d(t)] \tag{7a}$$

$$|v(t)| = \gamma(|F_d(t)| - \mu_s N_f) \tag{7b}$$

Using equation (5) we eliminate $x(t)$ in equations (4a) – (4d), and obtain the following equations:

$$\dot{v}(t) = \frac{1}{m}[F(t) - F_s(t) - F_d(t)] \tag{8a}$$

$$\dot{F}_s(t) = k\, v(t) \tag{8b}$$

$$0 = \beta(v(t), F_d(t)) \tag{8c}$$

The above system of equations is a system of *semi-explicit* differential-algebraic equations (DAEs) (for example, see Reference [2]). The first two equations are ordinary differential equations (ODEs), and the third equation is an algebraic constraint. It is, in general, not possible to find analytical solutions for these system of equations, and one must resort to numerical solutions. In a subsequent section, we shall present a numerical algorithm to solve these differential-algebraic equations. Some representative references on differential-algebraic equations are [2, 3, 6, 7, 8, 9].

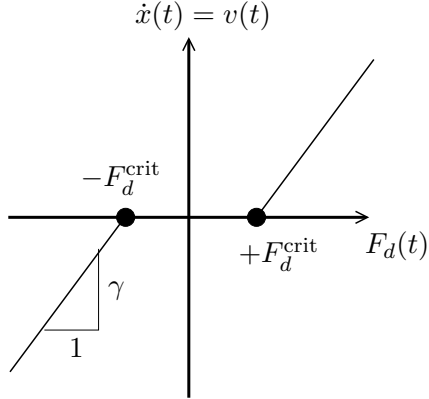

Figure 4. A pictorial description of Bingham model

**Remark 2.1.** *When $\gamma \to \infty$, Figure 4 resembles the response for Coulomb damping. However, it should be borne in mind that one cannot obtain the frictional response of solids by a limiting process for the frictional response of fluids. While the response curves might look the same, the philosophical underpinnings as well as other relevant physical issues are quite different.*

2.2. **Computation of a solution for the differential-algebraic system of equations.** Classical textbooks on vibrations such as the one written by Meirovitch supposedly provide a method for obtaining the response of the spring-mass system subject to Coulomb damping (see Meirovitch [11, pages 31-34]). The method consists of piecing together solutions corresponding to the regimes when $\dot{x} > 0$ and $\dot{x} < 0$. The resulting solution obtained by Meirovitch is smooth, as can be seen from the figure on page 34 of [11]. The procedure adopted by Meirovitch runs into trouble when one considers the following initial conditions: $x(0) = x_0$, where $\dot{x}(0) = 0$ and $|x_0| < \frac{\mu_s N_f}{k}$. The main difficulty lies in the inability to specify $F_d$ when $\dot{x}(0) = 0$.

Put differently, the main error of Meirovitch's analysis stems from the fact that when the initial conditions are such that $|x_0| < \mu_s N_f/k$ and $\dot{x}(0) = 0$, one cannot use the differential equation (1) to describe the system as $F_d$ is undefined, within the framework that Meirovitch used to describe the problem. The differential equation (1) is incapable of describing the response of the system. One can proceed in either one of two ways. Use the differential-algebraic system (1), (5) and (6). This framework does allow one to describe the system as we now do not choose to provide a constitutive relation for $F_d$, which we are not in a position to provide, rather we choose to treat the force $F_d$ as an unknown that is obtained as a solution to the system (1), (5) and (6). The other way of studying the problem that is well-posed is to describe the physical problem within the context of a differential inclusion, as considered by Filippov [5]. The approach adopted by Filippov [5] is described briefly in the discussion that follows. We shall, however, consider the problem within the context of a differential-algebraic system in some detail as such an approach can be adopted to study a whole host of problems wherein the forces in the springs and the dashpots are given in terms of implicit constitutive relations that relate the forces to appropriate kinematical variables.

If one is to interpret the frictional force in terms of the signum function, with $\text{sgn}[0] = 0$, then one does not have a classical solution to the problem. One may appeal to intuition in this case and claim that the solution is $x(t) = x_0$ and $\dot{x}(t) = 0$. However, this solution, as we will see, is not a solution in the classical sense. Neither does Meirovitch specify the sense in which the function he has constructed actually solves the problem. The absence of the specification of the sense to which one has a solution to the problem, namely the class of functions to which the solution belongs, can lend itself to misinterpretation of the function being a solution in the classical sense. Meirovitch takes an intuitive engineering approach to the problem and does not discuss the mathematical subtleties of the solution.

2.2.1. *Non-existence of classical solutions.* Before we proceed further, it is better to define what we mean by a classical solution or a solution in the sense of Caratheodory [4]. Let $z(t) \in \Re^n$ represent the state of a dynamical system at time $t$ and $f(t, z(t))$ denote a function mapping $\Re \times \Re^n \to \Re^n$. Let $z_0$ represent the state of the dynamical system at $t = 0$.

**Definition 2.2.** *Consider an $n^{th}$ degree ordinary differential equation with the following representation:*

(9) $$\dot{z}(t) = f(t, z(t)), \quad z(0) = z_0.$$

*We will say that $\phi(t)$ is a solution in the sense of Caratheodory to the differential equation (9) if*

- *$\phi(0) = z_0$,*
- *$\phi(t)$ is absolutely continuous, and*
- *$\phi(t) = z_0 + \int_0^t f(\tau, \phi(\tau)) d\tau$.*

*We will define $\phi(t)$ to be an equilibrium solution if, in addition to it being a solution, $\phi(t) = \phi(0)$ for all $t \geq 0$.*

The non-existence of classical solutions stems from the nature of dry friction. Physically, there is a continuum of equilibria which correspond to the velocity of the mass being equal to zero and the deflection of the spring not being sufficiently large for the spring force to overcome the friction force and cause motion. At all these equilibria, the acceleration and velocity of the mass is zero. However, if we were to set the acceleration to be zero when the velocity is zero, we see that the friction force must balance the spring force and hence, the friction force cannot be specified independently of the initial deflection of the spring. In fact, depending on the initial deflection of the spring, the Coulomb damping force can take any value in $(-\mu_s N_f, \mu_s N_f)$ when the velocity is zero. Hence, the specification of Coulomb damping force, $F_d$, cannot be performed solely on the basis of velocity as it is dependent on the deflection of the spring when $\dot{x} = 0$.

If one were to describe the damping force to be $F_d = -\mu_s N_f \, \text{sgn}[\dot{x}]$, the implication is that the damping force is specified when the velocity is zero and this is not physically the case. It is exactly this difficulty that leads to the non-existence of classical solutions for almost all initial conditions, where the initial velocity is zero and the initial deflection of the spring is no more than $\frac{\mu_s N_f}{k}$ in magnitude.

2.2.2. *Existence of solutions in the sense of Filippov.* Let us begin with the definition of a solution in the sense of Filippov [14]. We will use the following notation: For every $z \in \Re^n$, the term $B(z, \delta)$ represents the set of all points in $\Re^n$ which are no farther than $\delta$ from $z$. For every set $N \subset \Re^n$, let $\mu(N)$ denote the measure of the set. Let $A, B$ be two subsets of $\Re^n$; their difference, represented by $A - B$ indicates the set of points in $A$ that do not belong to $B$.

In Reference [5], Filippov considered the following differential inclusion:

$$\dot{x} \in F(x, t) \tag{10}$$

where $F(x, t)$ is the possible set of values that $\dot{x}$ can take. Filippov indicated that if for every $(x, t)$, $F(x, t)$ is non-empty and convex, and the distance between the sets $F(x', t')$ and $F(x, t)$ tends to zero as $x \to x'$ and $t \to t'$, then a solution (that is absolutely continuous and satisfies the differential inclusion (10) almost everywhere in time) exists.

Consider the set of values the right-hand side of equations (8a) and (8b) can take, which can be written as

$$F_+ = \begin{pmatrix} \frac{1}{m}\left[F(t) - F_s - (\frac{v}{\gamma} + \mu_s N_f)\right] \\ kv \end{pmatrix}, \text{ for } v(t) > 0 \tag{11a}$$

$$F_- = \begin{pmatrix} \frac{1}{m}\left[F(t) - F_s - (\frac{v}{\gamma} - \mu_s N_f)\right] \\ kv \end{pmatrix}, \text{ for } v(t) < 0 \tag{11b}$$

$$F_0 = \begin{pmatrix} \frac{1}{m}[F(t) - F_s - F_d] \\ 0 \end{pmatrix}, \text{ for } v(t) = 0, \text{ where } |F_d| \leq \mu_s N_f \tag{11c}$$

For the case $v \neq 0$, the set $F(x, t)$ is a singleton set. For the case $v = 0$, $F(x, t)$ is a convex set and the other conditions are met. Solutions exist and are guaranteed to pick appropriate value of $F_d$ from the set $\{|F_d| \leq \mu_s N_f\}$ when the velocity is zero.

2.3. **A numerical algorithm for solving a system of differential-algebraic equations.** It is well-known that Backward Difference Formulae (BDF) are stable and accurate (and hence more suitable) for solving DAEs numerically [2, 3]. The backward Euler scheme is the simplest member of BDF. We shall now discretize equations (8) using the backward Euler time stepping scheme. To this end, we shall discretize the time interval of interest into $N$ time instants denoted by $t_n$ $(n = 0, \cdots, N)$. For simplicity, we shall assume uniform time steps, and shall denote it by $\Delta t := t_n - t_{n-1}$. We shall denote the (time) discretized version of a given quantity $z(t)$ at the instant of time $t_n$ as

$$z^{(n)} \approx z(t = t_n) \quad n = 0, \cdots, N \tag{12}$$

The corresponding discretized equations using the backward Euler scheme at the instant of time $t_{n+1}$ can be written as

$$\frac{1}{\Delta t}[v^{(n+1)} - v^{(n)}] = \frac{1}{m}[F^{(n+1)} - F_s^{(n+1)} - F_d^{(n+1)}] \tag{13a}$$

$$\frac{1}{\Delta t}[F_s^{(n+1)} - F_s^{(n)}] = k\, v^{(n+1)} \tag{13b}$$

$$0 = \beta(v^{(n+1)}, F_d^{(n+1)}) \tag{13c}$$

Using equation (13b) we can write $F_s^{(n+1)}$ in terms of the variable $v^{(n+1)}$ as

$$F_s^{(n+1)} = F_s^{(n)} + \Delta t\, k\, v^{(n+1)} \tag{14}$$

Using the above equation, equation (13a) can be written as

$$\left[1 + \Delta t^2 \frac{k}{m}\right] v^{(n+1)} = \frac{\Delta t}{m}\left(\tilde{F}^{(n+1)} - F_d^{(n+1)}\right) \tag{15a}$$

$$\tilde{F}^{(n+1)} := \frac{m}{\Delta t} v_n + F^{(n+1)} - F_s^{(n)} \tag{15b}$$

From equation (15a) (and noting equation (7a) for the case of $|F_d(t)| > \mu_s N_f$) we have the following useful result

$$\text{sgn}[\tilde{F}^{(n+1)}] = \text{sgn}[F_d^{(n+1)}] \tag{16}$$

which is valid for all cases (that is, for both $|F_d^{(n+1)}| \le \mu_s N_f$ and $|F_d^{(n+1)}| > \mu_s N_f$). The above result can be utilized quite effectively in a numerical scheme as $\tilde{F}^{(n+1)}$ is known whereas $F_d^{(n+1)}$ is not known as it is part of the solution for the instant of time $t_{n+1}$. We shall employ a predictor-corrector-type scheme to calculate the values for the instant of time $t_{n+1}$ assuming that all the values at the prior instants of time are known. The scheme is outlined in Algorithm 1.

**Remark 2.3.** *In the case of the standard Bingham model, the (discrete) solution from $t_n$ to $t_{n+1}$ can be solved analytically. But for more complicated models, one may have to employ a Newton-Raphson scheme to solve the nonlinear equation, which can be written just in terms of $F_d^{(n+1)}$.*

**Algorithm 1** Predictor-corrector algorithm for advancing the solution from $t_n$ to $t_{n+1}$

---
1: Input: $k$, $m$, $\Delta t$, $\mu_s N_f$, $v^{(n)}$, $F_s^{(n)}$, $F^{(n+1)}$
2: Output: $x^{(n+1)}$, $v^{(n+1)}$, $F_d^{(n+1)}$, $F_s^{(n+1)}$
3: Predictor step: Calculate the predictor $\tilde{F}^{(n+1)} := \frac{m}{\Delta t}v^{(n)} + F^{(n+1)} - F_s^{(n)}$
4: Corrector step:
5: **if** $|\tilde{F}^{(n+1)}| \leq \mu_s N_f$ **then**
6: $\quad v^{(n+1)} = 0$ and $F_d^{(n+1)} = \tilde{F}^{(n+1)}$
7: **else**
8: $\quad F_d^{(n+1)} = \frac{1}{\gamma[1+\Delta t^2 \frac{k}{m}] + \frac{\Delta t}{m}} \left[ \gamma \left(1 + \Delta t^2 \frac{k}{m}\right) \operatorname{sgn}[\tilde{F}^{(n+1)}]\mu_s N_f + \frac{\Delta t}{m} \tilde{F}^{(n+1)} \right]$
9: $\quad v^{(n+1)} = \gamma(F_d^{(n+1)} - \operatorname{sgn}[F_d^{(n+1)}]\mu_s N_f)$
10: **end if**
11: Calculate $F_s^{(n+1)}$ and $x^{(n+1)}$ using
$$F_s^{(n+1)} = F_s^{(n)} + \Delta t\, k\, v^{(n+1)}, \quad x^{(n+1)} = \frac{1}{k} F_s^{(n+1)}$$

---

## 3. NUMERICAL EXAMPLES

In all the numerical examples we shall take $\gamma = 1$, $\mu_s N_f = 1$, $m = 1$, $k = 100$, and $\Delta t = 10^{-4}$.

### 3.1. Non-zero external forcing function.
We shall take $x(0) = 0$ and $F_d(0) = 0$. These conditions will imply that the corresponding (consistent) initial conditions for the other variables will be $v(0) = 0$ (from equation (6)) and $F_s(0) = 0$ (from equation (5)). Two different external forcing functions are considered

(17a) $$F_1(t) = \begin{cases} 0.5\sin(5\pi t) & 0 \leq t \leq 1, \\ 0 & t > 1, \end{cases}$$

(17b) $$F_2(t) = \begin{cases} 10.0\sin(5\pi t) & 0, \leq t \leq 1 \\ 0 & t > 1. \end{cases}$$

The numerical results for the forcing function $F_1(t)$ are presented in Figure 5. The numerical results for the forcing function $F_2(t)$ are presented in Figures 6 – 9.

**Remark 3.1.** *The forcing functions are chosen in such a way that*

(18) $$\max[F_1(t)] < \mu_s N_f \quad \text{and} \quad \max[F_2(t)] > \mu_s N_f$$

*The analytical solution for the first case is*

(19) $$F_d(t) = F_1(t),\ F_s(t) = 0,\ x(t) = 0,\ v(t) = 0$$

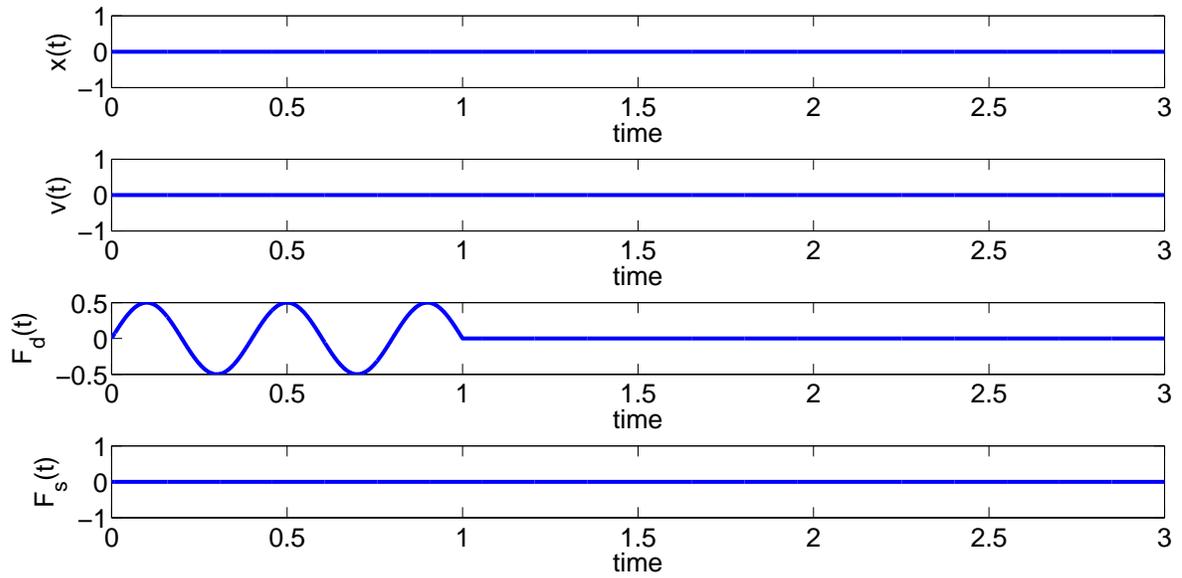

FIGURE 5. Non-zero forcing function using equation (17a).

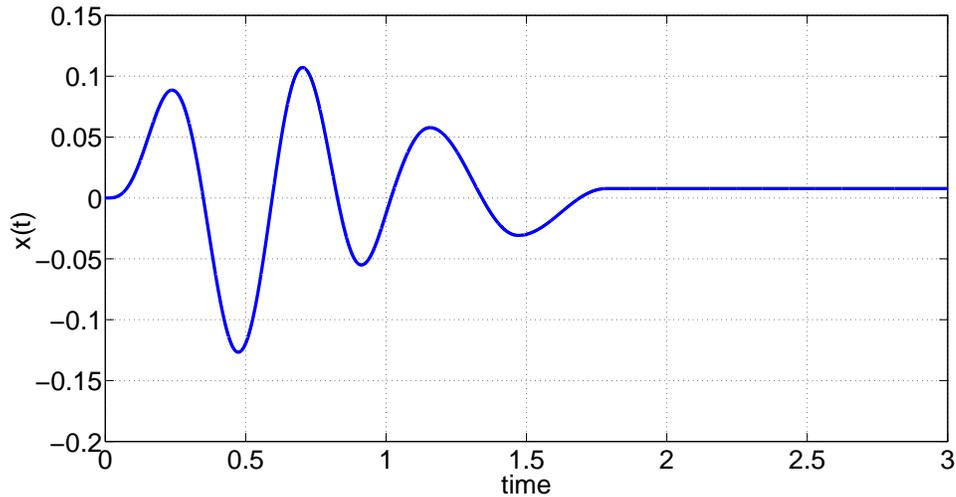

FIGURE 6. Displacement $x(t)$ for non-zero forcing function given in equation (17b).

*As one can see from Figure 17a, the numerical solution matches the analytical solution quite well.*

3.2. **Non-zero initial displacement.** We shall take the external force $F(t) = 0$, and take $F_d(0) = 0$, which implies that the consistent initial condition for the velocity will be $v(0) = 0$ (from equation (6)). Two different initial displacements are considered: $x(t = 0) = 0.005$ (which implies that $F_s(t = 0) = 0.5$), and $x(t = 0) = 0.5$ (which implies that $F_s(t = 0) = 50$). The numerical results corresponding to $x(0) = 0.005$ are presented in Figures 10. The numerical results for the case $x(0) = 0.5$ are presented in Figures 11 – 14.

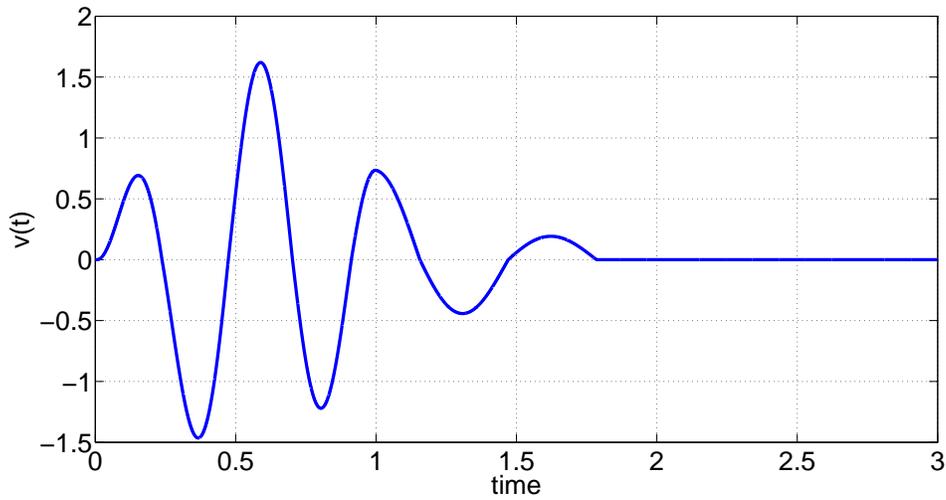

FIGURE 7. Velocity $v(t)$ using non-zero forcing function given in equation (17b).

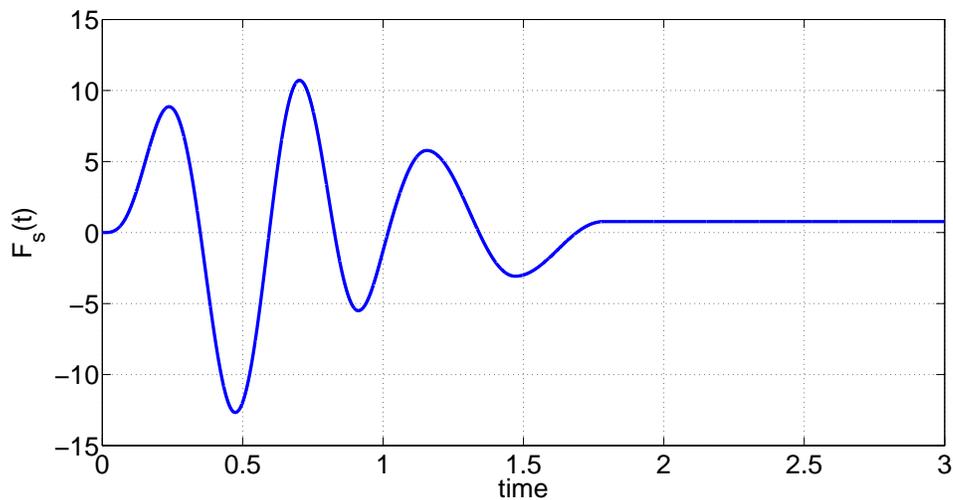

FIGURE 8. Force in the spring $F_s(t)$ using non-zero forcing function given in equation (17b).

## 4. CONCLUSIONS

In this short study we have considered problems of vibratory motion of lumped parameter systems wherein the forces in the springs and dashpots belonging to the system cannot be given as functions of the kinematical quantities; rather the kinematical quantities are expressed as functions of the forces. This leads to the system being described by a set of differential-algebraic equations. We illustrate the type of equations one has to deal with by considering a simple system, for the sake of simplicity of illustration. We considered a very simple mass-spring-dashpot configuration with the spring being a linear spring and a dashpot that contains a Bingham fluid. We subjected this system to a class of forcing functions and initial conditions and found the solutions to be in keeping with physical expectations. In the future, we plan to study the vibration of systems whose components are described by implicit constitutive relations between the forces and the kinematical variables which cannot be simplified either to the forces being given in terms of the kinematics or vice-versa.

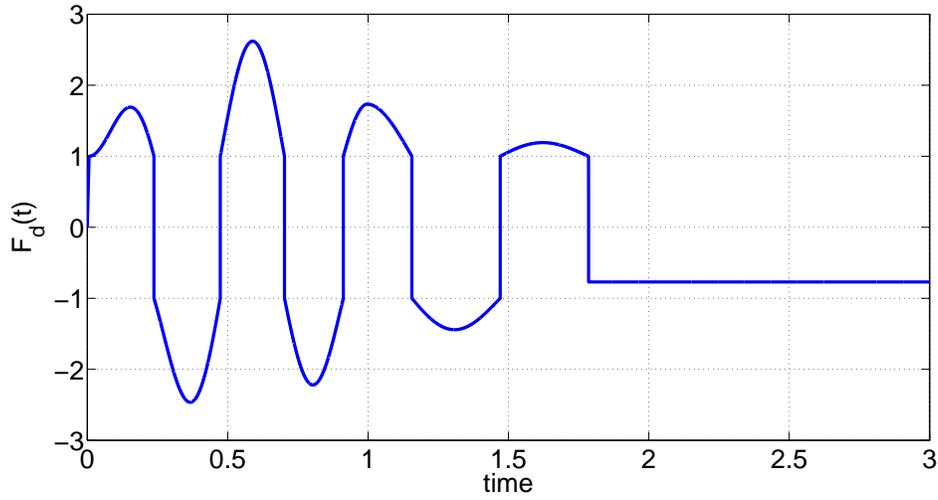

FIGURE 9. Force in the dashpot $F_d(t)$ using non-zero forcing function given in equation (17b).

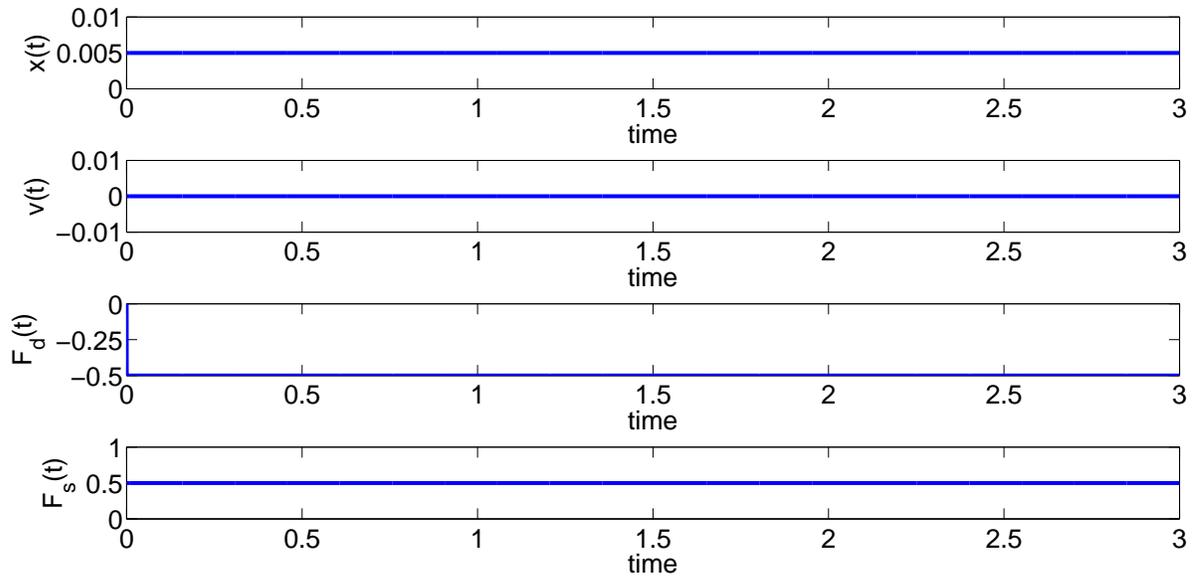

FIGURE 10. Non-zero initial displacement of $x(0) = 0.005$.

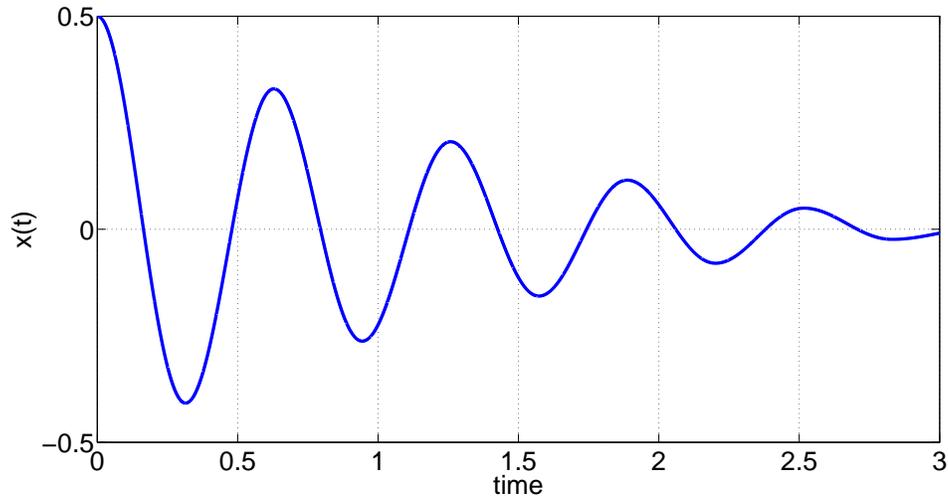

FIGURE 11. Displacement $x(t)$ using non-zero initial displacement of $x(0) = 0.5$.

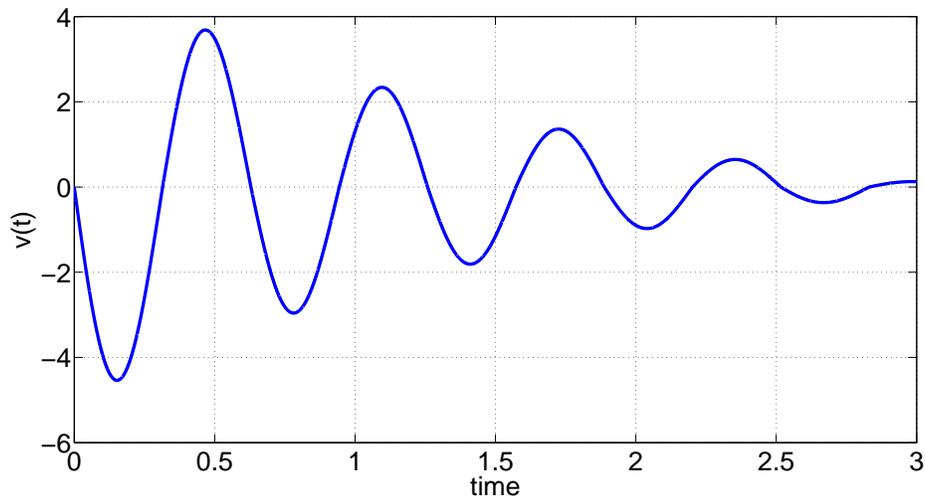

FIGURE 12. Velocity $v(t)$ using non-zero initial displacement of $x(0) = 0.5$.

## References

[1] E. Pa. Antonyuk. On models of dynamic systems with dry friction. *International Applied Mechanics*, 43:554–559, 2007.

[2] U. M. Ascher and L. R. Petzold. *Computer Methods for Ordinary Differential Equations and Differential-Algebraic Equations*. SIAM, Philadelphia, USA, 1998.

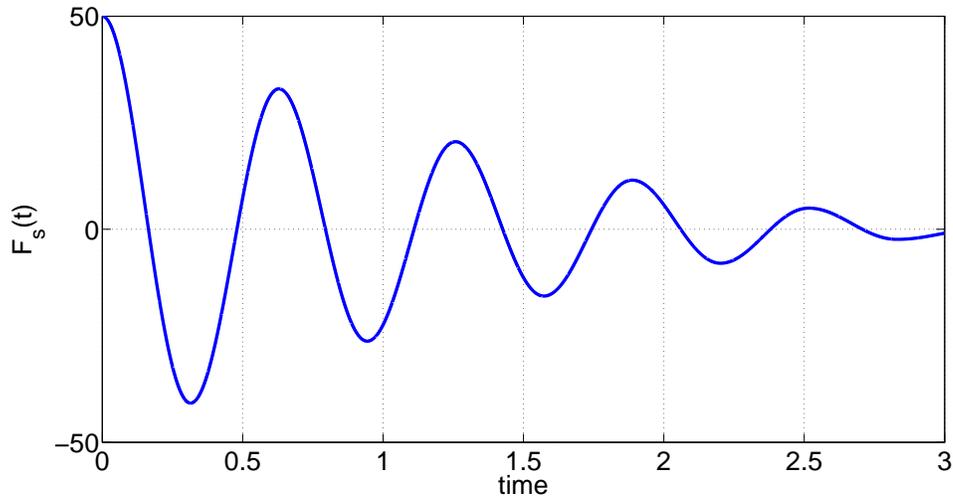

FIGURE 13. Force in the spring $F_s(t)$ using non-zero initial displacement of $x(0) = 0.5$.

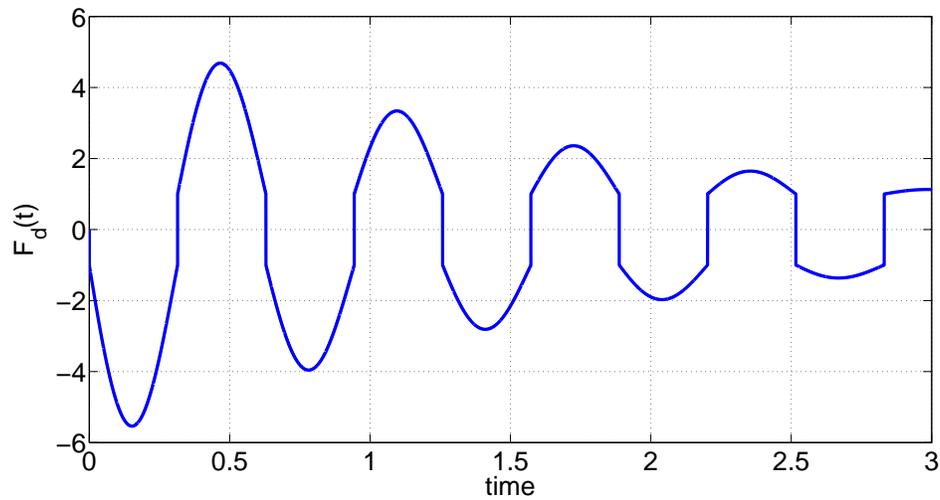

FIGURE 14. Force in the dashpot $F_d(t)$ using non-zero initial displacement of $x(0) = 0.5$.

[3] K. Brenan, S. Campbell, and L. Petzold. *Numerical Solutions of Initial-Value Problems in Differential-Algebraic Equations*. North-Holland, New York, USA, 1989.

[4] E. A. Coddington and N. Levinson. *Theory of Ordinary Differential Equations*. Krieger Publishing Company, Malabar, Florida, USA, 1984.

[5] A. F. Filippov. Classical solutions of differential equations with multivalued right-hand side. *SIAM Journal on Controls and Optimization*, 5:609–621, 1967.

[6] E. Hairer, C. Lubich, and M. Roche. *The Numerical Solution of Differential-Algebraic Systems by Runge-Kutta Methods*. Lecture Notes in Mathematics. Springer-Verlag, New York, USA, 1989.

[7] E. Hairer and G. Wanner. *Solving Ordinary Differential Equations II: Stiff and Differential-Algebraic Problems.* Springer-Verlag, New York, USA, 1996.

[8] P. Kunkel and V. Mehrmann. *Differential-Algebraic Equations.* European Mathematical Society, Zurich, Switzerland, 2006.

[9] R. M. M. Mattheij and J. Molenaar. *Ordinary Differential Equations in Theory and Practice.* John Wiley & Sons, Inc., Chichester, UK, 1997.

[10] J. C. Maxwell. On the dynamical theory of gases. *Philosophical Transactions of Royal Society of London*, A157:26–78, 1866.

[11] L. Meirovitch. *Elements of Vibration Analysis.* McGraw-Hill, Inc., New York, USA, 1986.

[12] K. R. Rajagopal. A generalized framework for studying the vibration of lumped parameter systems. *Submitted for Publication to Mechanics Research Communications.*

[13] K. R. Rajagopal. On implicit constitutive theories. *Application of Mathematics*, 28:279–319, 2003.

[14] S. Sastry. *Nonlinear Systems.* Springer-Verlag, New York, USA, 1999.

[15] M. Weircigroch and P. A. Zhilin. On the Painleve Paradoxes: Nonlinear Oscillations in Mechanical Systems. *Proceedings of the XXVII Summer Schools, St. Petersburg*, pages 1–22, 2007.

Dr. Swaroop Darbha, Department of Mechanical Engineering, 320 Engineering/Physics Building, Texas A&M University, College Station, Texas 77843. TEL:+1-979-862-2238

*E-mail address*: dswaroop@tamu.edu

Dr. Kalyana Babu Nakshatrala, Department of Mechanical Engineering, 216 Engineering/Physics Building, Texas A&M University, College Station, Texas 77843. TEL:+1-979-845-1292

*E-mail address*: knakshatrala@tamu.edu

Correspondence to: Dr. K. R. Rajagopal, Department of Mechanical Engineering, 314 Engineering/Physics Building, Texas A&M University, College Station, Texas 77843. TEL:+1-979-862-4552

*E-mail address*: krajagopal@tamu.edu